\documentclass{amsart}

\usepackage{amsmath}
\usepackage{amscd}
\usepackage{amssymb}

\newcommand{\bC}{{\mathbb C}}

\newcommand{\cS}{{\mathcal S}}
\newcommand{\Mbar}{\overline{\mathcal M}}
\newcommand{\cC}{{\mathcal C}}

\DeclareMathOperator{\Aut}{Aut}
\DeclareMathOperator{\id}{id}
\DeclareMathOperator{\tr}{tr}

\newtheorem{theorem}{Theorem}[section]
\newtheorem{theorem/definition}{Theorem/Definition}[section]
\newtheorem{proposition}{Proposition}[section]
\newtheorem{lemma}{Lemma}[section]

\theoremstyle{remark}
\newtheorem{remark}{Remark}[section]

\theoremstyle{definition}
 \newtheorem{example}{Example}[section]

\begin{document}

\title{Hodge integrals, Hurwitz numbers, and Symmetric Groups}
\author{Jian Zhou}
\address{Department of Mathematical Sciences\\Tsinghua University\\Beijing, 100084,
China}
\email{jzhou@math.tsinghua.edu.cn}
\footnote{This research is partially supported by research grants from CNSF and Tsinghua University.}
\begin{abstract}
We prove some combinatorial results related to a formula on Hodge integrals conjectured by Mari\~no and
Vafa.
These results play important roles in the proof and applications of this formula by the author jointly
with Chiu-Chu Melissa Liu and Kefeng Liu.
We also compare with some related results on Hurwitz numbers and obtain
some closed expressions for the generating series of Hurwitz numbers and the related Hodge integrals.
\end{abstract}

\maketitle

\section{Introduction}

In this paper we study some combinatorial problems related to the
conjecture of Mari\~{n}o-Vafa on Hodge integrals,
recently proved in joint work with Liu and Liu \cite{LLZ1, LLZ2}.
Some of the results proved here are used in the proof.
We also study the closely related Hurwitz numbers.
In both cases,
one has formulas expressing some special types of Hodge integrals
in terms of the characters of the symmetric groups.
The latter are closely related to other branches of mathematics,
e.g.,
symmetric functions,
representations of Lie algebras of type $A$ etc.
Hence understanding the relationship between the Hodge integrals and symmetric groups is an important step
towards establishing connections of the moduli spaces of Riemann surfaces with these mathematical branches.
Some of the results have been announced in \cite{Zho}.
See also \cite{LLZ1}.

 Hodge integrals are integrals on Deligne-Mumford
moduli space of curves of the following form:
$$\int_{\Mbar_{g,n}} \prod_{i=1}^n \psi_i^{m_i} \prod_{j=1}^g \lambda_j^{n_j}.$$
Their explicit evaluations are difficult problems important to both algebraic geometry and mathematical physics.
It is well known that the Hodge integrals can be reduced to integrals of the form
$$\int_{\Mbar_{g,n}} \prod_{i=1}^n \psi_i^{m_i}$$
which appear in two dimensional quantum gravity,
whose evaluations are given by the famous Witten's conjecture/Kontsevich's theorem \cite{Wit, Kon}.
See \cite{Fab} for the description of a {\em Maple} program  which in principle can be used to compute any Hodge integral.

In this work we will focus on some Hodge integrals that
appear in the calculations of Gromov-Witten invariants by localization techniques
\cite{Kon2, Gra-Pan} in the positive genus case.
More precisely,
we will be concerned with the following Hodge integrals appearing in the formal calculations of
some open Gromov-Witten invariants by localization method \cite{Kat-Liu}:
\begin{equation} \label{eqn:Kat-Liu}
\begin{split}
& (\tau(\tau+1))^{l(\mu)-1}\cdot \prod_{i=1}^{l(\mu)}
\frac{\prod_{j=1}^{\mu_i-1}(j+\mu_i\tau)}{(\mu_i-1)!} \cdot
\int_{\Mbar_{g, l(\mu)}} \frac{\Lambda^{\vee}_g(1)\Lambda^{\vee}_g(-\tau-1)\Lambda_g^{\vee}(\tau)}
{\prod_{i=1}^{l(\mu)} (1 - \mu_i\psi_i)}.
\end{split}
\end{equation}
We will refer to such Hodge integrals as conifold Hodge integrals.
In \cite{Mar-Vaf}, Mari\~{n}o and Vafa derived a closed
formula for a generating function of certain open Gromov-Witten invariants
by duality with Chern-Simon theory.
This leads to a conjectural formula of Hodge integrals by comparing with calculations in \cite{Kat-Liu}.
See \S \ref{sec:MV} for the precise formulation.
Note the two sides of the conjectured formula of Mari\~{n}o and  Vafa
are both mathematically well understood.
As mentioned above,
this formula has recently been proved in joint work with Liu and Liu \cite{LLZ1, LLZ2}.
In this work we will prove some combinatorial results that will be used in that proof.

A key result in this work
is that the right-hand side of the Mari\~no-Vafa formula satisfies the cut-and-join equations.
This is motivated by similar results for Hurwitz numbers,
for which there are two previous proofs: by combinatorial method \cite{Gou-Jac-Vai}
and by geometric method \cite{Li-Zha-Zhe, Ion-Par}.
Combinatorial method is used in this paper to prove the cut-and-join equation for the combinatorial expressions
that appear on the right-hand side of  the Mari\~no-Vafa formula,
and geometric method is used in \cite{LLZ1, LLZ2} to show the conifold Hodge integrals
on the left-hand side of the formula satisfy the same equation.
The exact match between the geometry of the moduli space and the combinatorics involved
is a salient feature of the proof presented in \cite{LLZ1, LLZ2}.
The ELSV formula \cite{ELSV, Gra-Vak} shows that the Hurwitz numbers are related to Hodge integrals
similar to the conifold Hodge integrals.
See \cite{LLZ3} for a proof of this formula along the same lines.

To make the analogies between the Hurwitz numbers and the conifold Hodge integrals more close,
we derive a formula similar to the Mari\~{n}o-Vafa formula for the Hurwitz numbers
by the Burnside formula.
It is clear that Mari\~no-Vafa formula is much more general than the ELSV formula,
and it by no means follows easily from the latter.
As remarked by Mari\~no and Vafa,
it is in principle possible to derive all Hodge integrals
containing up to three $\lambda$ classes from this formula,
but it is not clear how to do so from the latter.
See \cite{LLZ3} for some examples.
We also show how to obtain closed formulas for generating series of Hurwitz numbers
by the cut-and-join equations.
This not only serves as the prototype for the proof in \cite{LLZ1, LLZ2},
but also has its own independent interest.

This paper is a revised version of a preprint which was limitedly circulated
since November 2002.
An earlier version is available from my homepage since June 2003.

\section{Preliminaries}
\label{sec:Pre}

In this section we recall some basic facts from the combinatorics of symmetric groups,
including a technical result used to establish the cut-and-join equations.

\subsection{Partitions}

A partition of a positive integer $d$ is a sequence of integers $\mu_1 \geq \mu_2 \geq \dots \geq \mu_l > 0$
such that
$$\mu_1 + \dots + \mu_l = d.$$
We write
\begin{align*}
|\mu| & = d, &
l(\mu) & = l.
\end{align*}
Denote by $m_j(\mu)$ the number of $j$'s among $\mu_1, \dots, \mu_l$.
The automorphism group $\Aut(\mu)$ of $\mu$ consists of possible permutations among the $\mu_i$'s,
hence its order is given by
$$|\Aut(\mu)| = \prod_j m_j(\mu)!.$$
We will also need the following number associated to a partition:
\begin{eqnarray} \label{eqn:n1}
&& n(\mu) = \sum_i (i-1)\mu_i.
\end{eqnarray}
The following number for a partition will be used in the Mari\~{n}o-Vafa formula:
\begin{eqnarray*}
&& \kappa_{\mu} = \sum_{i=1}^{l(\mu)} \mu_i(\mu_i-2i+1).
\end{eqnarray*}

Each partition $\mu$ of $d$ corresponds to a conjugacy class $C(\mu)$ of the symmetric group $S_d$.
For simplicity of notations,
we denote by $C(2)$ the conjugacy class of transpositions.
The number of elements in $C(\mu)$ is given by
$$|C(\mu)| = \frac{d!}{z_{\mu}},$$
where
$$z_{\mu} = \prod_j m_j(\mu)!j^{m_j(\mu)} .$$
Denote by $(-1)^g$ the sign of an element in $S_d$.
It is easy to see that
$$(-1)^g = (-1)^{|\mu| - l(\mu)},$$
for $g \in C(\mu)$.

The Young diagram of $\mu$ has $m_j(\mu)$ rows of squares of length $j$.
The partition corresponding to the transpose of the Young diagram of $\mu$
will be denoted by $\mu'$.
Denote by $h(e)$ the hook length of a square $e \in \mu$.
We will need the following identity (\cite[(1.6)]{Mac}):
\begin{eqnarray}
&& n(\mu) = \sum_i \begin{pmatrix} \mu'_i \\2 \end{pmatrix}.  \label{eqn:n2}
\end{eqnarray}

\subsection{Representations of symmetric groups}
Each partition $\nu$ corresponds to an irreducible representation $R_{\nu}$ of $S_d$.
The value of the character $\chi_{R_{\nu}}$ on the conjugacy class $C(\mu)$ is
denoted by $\chi_{\nu}(C(\mu))$.
They can be calculated using symmetric functions by the following formula \cite{Mac}:
$$s_{\nu} = \sum_{|\mu|=|\nu|} \frac{\chi_{\nu}(C(\mu))}{z_{\mu}} p_{\mu}.$$
For example,
$\chi_{(d)}$ corresponds to the trivial representation,
$\chi_{(1^d)}$ corresponds to the sign representation.

Set
$$c_{\mu} = \sum_{g \in C(\mu)} g.$$
Since $c_{\mu}$ lies in the center of the group algebra $\bC S_d$,
it acts as an scalar $f_{\nu}(\mu)$ on any irreducible representation $R_{\nu}$.
Now
$$\dim R_{\nu} \cdot f_{\nu}(\mu) = \tr|_{R_{\nu}} c_{\mu}
= \sum_{g \in C(\mu)} \tr|_{R_{\nu}} g
= |C(\mu)| \chi_{\nu}(\mu),$$
hence
$$f_{\nu}(\mu) = |C(\mu)|\frac{\chi_{\nu}(\mu)}{\dim R_{\nu}}
= \frac{d!}{\dim R_{\nu}} \cdot  \frac{\chi_{\nu}(C(\mu))}{z_{\mu}}.$$
We summarize the properties of these numbers in the following:

\begin{proposition} \label{prop:f}
Each $f_{\nu}(\mu)$ is an integer.
Furthermore, one has
\begin{eqnarray}
&& f_{\nu^t}(\mu) = (-1)^{|\mu|-l(\mu)}f_{\nu}(\mu),  \label{eqn:trans} \\
&& f_{\nu}(2) = \frac{1}{2} \kappa_{\nu}. \label{eqn:kappa}
\end{eqnarray}
\end{proposition}

\begin{proof}
For a proof of the first statement,
see e.g. \cite{Mac}, p. 126, Example 16,
where $f_{\nu}(\mu)$ is written as $\omega^{\nu}_{\mu}$.
For (\ref{eqn:trans}), recall
$$\chi_{\nu^t} = \chi_{\nu}\chi_{(d)},$$
hence
\begin{eqnarray*}
&& \chi_{\nu^t}(\mu) = (-1)^{|\mu|-l(\mu)}\chi_{\nu}(\mu).
\end{eqnarray*}
See e.g. \cite{Mac}, p. 116, Example 2.
Hence we have
\begin{eqnarray*}
&& \dim R_{\nu^t} = \chi_{\nu^t}(1^{|\nu^t|}) = \chi_{\nu}(1^{|\nu|}) = \dim R_{\nu}, \\
&& f_{\nu^t}(\mu) = |C(\mu)|\frac{\chi_{\nu^t}(\mu)}{\dim R_{\nu^t}}
= (-1)^{|\mu|-l(\mu)}|C(\mu)|\frac{\chi_{\nu}(\mu)}{\dim R_{\nu}}
= (-1)^{|\mu|-l(\mu)}f_{\nu}(\mu).
\end{eqnarray*}
By \cite{Mac}, p. 118, Example 7,
\begin{eqnarray*}
f_{\nu}(2) & = & |C(2)| \frac{\chi_{\nu}(C(2))}{\dim R_{\nu}} = n(\nu') - n(\nu)\\
& = & \sum_{i=1}^{l(\nu)} \begin{pmatrix} \nu_i \\ 2\end{pmatrix}
- \sum_{i=1}^{l(\nu)} (i-1)\nu_i \\
& = & \frac{1}{2} \sum_{i=1}^{l(\nu)} (\nu_i^2 - 2i \nu_i + \nu_i) = \frac{1}{2} \kappa_{\nu}.
\end{eqnarray*}
In the above we have used (\ref{eqn:n1}) and (\ref{eqn:n2}).
\end{proof}

\subsection{Cut-and-join analysis}
Suppose $\mu$ is a partition containing $i$ and $j$,
the partition $\nu$ obtained from $\mu$ by replacing $i,j$ with $i+j$ is called a $(i,j)$-join of $\mu$,
and $\mu$ is called an ($i, j)$-cut of $\nu$.
We will write $\nu \in  J_{i,j}(\mu)$ and $\mu \in C_{i, j}(\nu)$.
It is easy to see that
\begin{eqnarray} \label{eqn:CJm}
&& \frac{m_{i+j}(\nu)}{\prod_k m_k(\nu)!} = \begin{cases}
\frac{m_i(\mu)m_j(\mu)}{\prod_k m_k(\mu)!}, & i \neq j, \\
\frac{m_i(\mu)(m_i(\mu)-1)}{\prod_k m_k(\mu)!}, & i =j.
\end{cases}
\end{eqnarray}

Denote by $[s_1, \dots, s_k]$ a $k$-cycle.
Then
$$[s, t]\cdot [s, s_2, \dots, s_i, t, t_2, \dots, t_j] = [s, s_2, \dots, s_i][t, t_2, \dots, t_j],$$
i.e., an $i+j$-cycle is cut into an $i$-cycle and a $j$-cycle.
Conversely,
$$[s, t]\cdot [s, s_2, \dots, s_i][t, t_2, \dots, t_j] = [s, s_2, \dots, s_i, t, t_2, \dots, t_j],$$
i.e., an $i$-cycle and a $j$-cycle is joined to an $i+j$-cycle.
To summarize,
the cuts and joins can be realized by multiplications by transpositions.

\begin{lemma}
Suppose $h \in S_d$ has cycle type $\mu$.
The product $C_{(2)} \cdot h$ is a sum of elements of $S_d$ whose type is either a cut or a join of $\mu$.
More precisely,
there are $ijm_i(\mu)m_j(\mu)$ (when $i < j$) or
$i^2 m_i(\mu)(m_i(\mu)-1)/2$ (when $i=j$) elements obtained from $h$
by joining an $i$-cycle in $h$ to a $j$-cycle in $h$,
and there are $(i+j)m_{i+j}(\mu)$ (when $i < j$) or $im_{2i}(\mu)$ (when $i=j$) elements obtained from $h$ by
cutting an $(i+j)$-cycle into an $i$-cycle and a $j$-cycle.
\end{lemma}

\begin{proof}
For a permutation $h$ of type $\mu$,
$c_{(2)} \cdot h$ is a sum of all elements obtained from $h$ by either a cut or a join.
Fix a pair of $i$-cycle and $j$-cycle of $h$,
there are $i \cdot j$ different ways to join them to an $(i+j)$-cycle.
Taking into the account of $m_i(\mu)$ choices of $i$-cycles, and $m_j(\mu)$ choices of $j$-cycles,
we get the number of different ways to obtain an element from $h$ by joining an $i$-cycle in $h$ to a $j$-cycle in $h$
is
$$\begin{cases}
ijm_i(\mu)m_j(\mu), & i < j \\
i^2 m_i(\mu)(m_i(\mu)-1)/2, & i = j. \end{cases}$$
Similarly,
fix an $(i+j)$-cycle of $h$,
there are $i+j$ different ways to cut it into an $i$-cycle and a disjoint $j$-cycle in $h$.
And taking into account the number of $(i+j)$-cycles in $h$,
we get the number of different ways to obtain an element from $h$ by cutting an $(i+j)$-cycle
into an $i$-cycle and a $j$-cycle is
$$\begin{cases}
(i+j)m_{i+j}(\mu), & i < j, \\
im_{2i}(\mu), & i = j.\end{cases}$$
\end{proof}

The following result will be used to establish the cut-and-join equations.

\begin{proposition} \label{prop:CJ}
For any partition $\nu$ of $d$,
one has
\begin{eqnarray*}
&& f_{\nu}(2)  \cdot \sum_{\mu}  \frac{\chi_{\nu}(\mu)}{z_{\mu}}p_{\mu}
= \frac{1}{2} \sum_{i,j}
\left((i+j)p_ip_j \frac{\partial}{\partial p_{i+j}}
+ ijp_{i+j}\frac{\partial}{\partial p_i}\frac{\partial}{\partial p_j}\right)
\sum_{\eta} \frac{\chi_{\nu}(\eta)}{z_{\eta}} p_{\eta}.
\end{eqnarray*}
\end{proposition}

\begin{proof}
For any $h \in S_d$ of cycle type $\mu$ we have
\begin{eqnarray*}
&& \sum_{\mu}  f_{\nu}(2) \frac{\chi_{\nu}(\mu)}{z_{\mu}}p_{\mu} \\
& = & \sum_{\mu} \tr [f_{\nu}(2) \id \cdot \rho_{\nu}(h)]
\cdot  \prod_i \frac{p_i^{m_i(\mu)}}{i^{m_i(\mu)} m_i(\mu)!} \\
& = & \sum_{\mu} \tr [\sum_{g \in C(2)} \rho_{\nu}(g) \cdot \rho_{\nu}(h)]
\cdot  \prod_i \frac{p_i^{m_i(\mu)}}{i^{m_i(\mu)} m_i(\mu)!} \\
& = &  \sum_{\mu}\tr \rho_{\nu}(\sum_{g \in C(2)} g\cdot h)
\cdot  \prod_i \frac{p_i^{m_i(\mu)}}{i^{m_i(\mu)} m_i(\mu)!} \\
& = & \sum_{\mu}\left( \sum_{i < j} \left(\sum_{\eta \in J_{i,j}(\mu)} ijm_i(\mu)m_j(\mu) \chi_{\nu}(\eta)
+ \sum_{\eta \in C_{i,j}(\mu)} (i+j)m_{i+j}(\mu) \chi_{\nu}(\eta)\right) \right. \\
&& \left. + \sum_i\left(\sum_{\eta \in J_{i,i}(\mu)} \frac{1}{2}i^2m_i(\mu)(m_i(\mu)-1) \chi_{\nu}(\eta)
+ \sum_{\eta \in C_{i,i}(\mu)} im_{2i}(\mu) \chi_{\nu}(\eta)\right)\right) \\
&& \cdot  \prod_i \frac{p_i^{m_i(\mu)}}{i^{m_i(\mu)} m_i(\mu)!} \\
& = & \frac{1}{2} \sum_{i,j}
\left((i+j)p_ip_j \frac{\partial}{\partial p_{i+j}}
+ ijp_{i+j}\frac{\partial}{\partial p_i}\frac{\partial}{\partial p_j}\right)
\sum_{\eta} \frac{\chi_{\nu}(\eta)}{z_{\eta}} p_{\eta}.
\end{eqnarray*}
In the last equality we have used (\ref{eqn:CJm}).
\end{proof}

\begin{remark}
The differential operator in the preceding Proposition first appeared \cite{Gou},
proved by similar argument using the characteristic map.
See also \cite{Fre-Wan}.
\end{remark}

\section{Mari\~{n}o-Vafa Formula and Some Related Combinatorial Problems}

\label{sec:MV}

In this section we give the precise formulation of the Mari\~no-Vafa formula,
and prove some combinatorial results related to it.

\subsection{Mari\~{n}o-Vafa formula}

For a partition $\mu$  define
\begin{eqnarray*}
\cC_{g, \mu}(\tau)& = & - \frac{\sqrt{-1}^{|\mu|+l(\mu)}}{|\Aut(\mu)|}
 [\tau(\tau+1)]^{l(\mu)-1}
\prod_{i=1}^{l(\mu)}\frac{ \prod_{a=1}^{\mu_i-1} (\mu_i \tau+a)}{(\mu_i-1)!} \\
&& \cdot \int_{\Mbar_{g, l(\mu)}}
\frac{\Lambda^{\vee}_g(1)\Lambda^{\vee}_g(-\tau-1)\Lambda_g^{\vee}(\tau)}
{\prod_{i=1}^{l(\mu)}(1- \mu_i \psi_i)}, \\
\cC_{\mu}(\lambda; \tau) & = & \sum_{g \geq 0} \lambda^{2g-2+l(\mu)}\cC_{g,
\mu}(\tau)
\end{eqnarray*}

Note that
$$
\int_{\Mbar_{0, l(\mu)}}
\frac{\Lambda^\vee_0(1)\Lambda^\vee_0(-\tau-1)\Lambda_0^\vee(\tau) }
{\prod_{i=1}^{l(\mu)} (1 -\mu_i \psi_i) }
=\int_{\Mbar_{0, l(\mu)}}
\frac{1}{\prod_{i=1}^{l(\mu)}(1 - \mu_i\psi_i)}
= |\mu|^{l(\mu)-3}
$$
for $l(\mu)\geq 3$, and we use this expression to extend the definition
to the case $l(\mu)<3$.

Introduce formal variables $p=(p_1,p_2,\ldots,p_n,\ldots)$, and define
$$
p_\mu=p_{\mu_1}\cdots p_{\mu_{l(\mu)} }
$$
for a partition $\mu=(\mu_1\geq \cdots \geq \mu_{l(\mu)}>0 )$. Define
generating functions
\begin{eqnarray*}
\cC(\lambda; \tau; p) & = & \sum_{|\mu| \geq 1} \cC_{\mu}(\lambda;\tau)p_{\mu}, \\
\cC(\lambda; \tau; p)^{\bullet} & = & e^{\cC(\lambda; \tau; p)}.
\end{eqnarray*}

As pointed out in \cite{Mar-Vaf},
by comparing physical computations in \cite{Mar-Vaf} with localization computations in \cite{Kat-Liu},
one obtains a conjectural formula for $\cC_{\mu}(\tau)$.
This formula can be explicitly written down as follows:
\begin{equation}\label{eqn:Mar-Vaf1}
\cC(\lambda; \tau; p)
=\sum_{n \geq 1} \frac{(-1)^{n-1}}{n}\sum_{\mu}
\left(\sum_{\cup_{i=1}^n \mu^i = \mu}
\prod_{i=1}^n \sum_{|\nu^i|=|\mu^i|} \frac{\chi_{\nu^i}(C(\mu^i))}{z_{\mu^i}}
e^{\sqrt{-1}(\tau+\frac{1}{2})\kappa_{\nu^i}\lambda/2} V_{\nu^i}(\lambda)
\right)p_\mu,
\end{equation}
\begin{equation}\label{eqn:Mar-Vaf2}
 \cC(\lambda;\tau; p)^{\bullet} = \sum_{|\mu| \geq 0}
\left(\sum_{|\nu|=|\mu|} \frac{\chi_{\nu}(C(\mu))}{z_{\mu}}
e^{\sqrt{-1}(\tau+\frac{1}{2})\kappa_{\nu}\lambda/2} V_{\nu}(\lambda)\right)
p_\mu,
\end{equation}
where
\begin{equation} \label{eqn:V}
\begin{split}
V_{\nu}(\lambda) = & \prod_{1 \leq a < b \leq l(\nu)}
\frac{\sin \left[(\nu_a - \nu_b + b - a)\lambda/2\right]}{\sin \left[(b-a)\lambda/2\right]} \\
& \cdot\frac{1}{\prod_{i=1}^{l(\nu)}\prod_{v=1}^{\nu_i} 2 \sin \left[(v-i+l(\nu))\lambda/2\right]}.
\end{split} \end{equation}

\subsection{A simple expression for $V_{\nu}(\lambda)$}

By writing down some low degree examples we discover a simple expression which we now come to.
Let us recall some notations for partitions.
The hook length of $\nu$ at the square $x$ located at the $i$-th row and $j$-th column
is defined to be:
$$h(x) = \nu_i + \nu_j'-i-j+1.$$
Then one has the following  two identities (\cite{Mac}, p.p. 10 - 11):
\begin{eqnarray}
&& \prod_{x \in \nu} (1 - t^{h(x)})
= \frac{\prod_{i=1}^{l(\nu)} \prod_{j=1}^{\nu_i - i + l(\nu)}( 1 - t^j)}
{\prod_{i < j} ( 1- t^{\nu_i - \nu_j - i +j})}, \label{eqn:hook1} \\
&& \sum_{x \in \nu} h(x) = n(\nu) + n(\nu') + |\nu|. \label{eqn:hook2}
\end{eqnarray}
With these preparations we can now prove the following result announced
in \cite[Theorem 1]{Zho}:

\begin{theorem}
We have
\begin{eqnarray}
&& V_{\nu}(\lambda) = \frac{1}{2^l\prod_{x \in \nu} \sin [h(x)\lambda/2]}.
\end{eqnarray}
\end{theorem}

\begin{proof}
We begin by rewriting the right-hand side of (\ref{eqn:hook1}) as follows.
\begin{eqnarray*}
&& RHS \\
& = &
\frac{\prod_{i=1}^{l(\nu)}\prod_{j=1}^{- i + l(\nu)}( 1 - t^j)}
    {\prod_{i < j} ( 1- t^{\nu_i - \nu_j - i + j})}
\cdot \prod_{i=1}^{l(\nu)}\prod_{j=1-i+l(\nu)}^{\nu_i - i + l(\nu)}( 1 - t^j) \\
& = &
\frac{\prod_{i<j} ( 1 - t^{j-i})}{\prod_{i < j} ( 1- t^{\nu_i - \nu_j - i +j})}
\cdot \prod_{i=1}^{l(\nu)}\prod_{j=1}^{\nu_i}( 1 - t^{j-i+l(\nu)}) \\
& = & t^{\left(\sum_{i<j} (j-i) - \sum_{i<j}(\nu_i-\nu_j-i+j) + \sum_{i=1}^{l(\nu)} \sum_{j=1}^{\nu_i} (j-i+l(\nu))\right)/2} \\
&& \cdot
\frac{\prod_{i<j} ( t^{-(j-i)/2} - t^{(j-i)/2})}{\prod_{i < j} ( t^{-(\nu_i-\nu_j-i+j)/2}- t^{(\nu_i - \nu_j - i +j)/2})}
\prod_{i=1}^{l(\nu)}\prod_{j=1}^{\nu_i}( t^{-(j-i+l(\nu))/2} - t^{(j-i+l(\nu))/2}).
\end{eqnarray*}
Now
\begin{eqnarray*}
&& \sum_{i<j} (j-i) - \sum_{i<j}(\nu_i-\nu_j-i+j) + \sum_{i=1}^{l(\nu)} \sum_{j=1}^{\nu_i} (j-i+l(\nu)) \\
& = & - \sum_{i <j}\nu_i + \sum_{i< j} \nu_j + \sum_{i=1}^{l(\nu)} \sum_{j=1}^{\nu_i} j
- \sum_{i=1}^{l(\nu)} \sum_{j=1}^{\nu_i} i + \sum_{i=1}^{l(\nu)} \sum_{j=1}^{\nu_i} l(\nu) \\
& = & - \sum_{i=1}^{l(\nu)} (l(\nu) - i) \nu_i + \sum_{j=1}^{l(\nu)} (j-1)\nu_j
+ \sum_{i=1}^{l(\nu)} \frac{\nu_i(\nu_i+1)}{2} -\sum_{i=1}^{l(\nu)} i\nu_i + |\nu|l(\nu)\\
& = & - |\nu|l(\nu) + \sum_{i=1}^{l(\nu)} i\nu_i + \sum_{j=1}^{l(\nu)} (j-1)\nu_j
+ \sum_{i=1}^{l(\nu)} \frac{\nu_i(\nu_i-1)}{2} + |\nu| - \sum_{i=1}^{l(\nu)} i\nu_i + |\nu|l(\nu) \\
& = &  \sum_{j=1}^{l(\nu)} (j-1)\nu_j + \sum_{i=1}^{l(\nu)} \frac{\nu_i(\nu_i-1)}{2} + |\nu| \\
& = & n(\nu) + n(\nu') + |\nu| \\
& = & \sum_{x \in \nu} h(x).
\end{eqnarray*}
Comparing with the left-hand side,
one then gets:
\begin{eqnarray*}
&& \prod_{x \in \nu} (t^{-h(x)/2} - t^{h(x)/2}) \\
& = & \frac{\prod_{i<j} ( t^{-(j-i)/2} - t^{(j-i)/2})}{\prod_{i < j} ( t^{-(\nu_i-\nu_j-i+j)/2}- t^{(\nu_i - \nu_j - i +j)/2})}
\prod_{i=1}^{l(\nu)}\prod_{j=1}^{\nu_i}( t^{-(j-i+l(\nu))/2} - t^{(j-i+l(\nu))/2}).
\end{eqnarray*}
The proof is completed by taking $t = e^{-\sqrt{-1}\lambda}$.
\end{proof}

\subsection{The cut-and-join equation}

\begin{theorem}
Denote by the right-hand sides of (\ref{eqn:Mar-Vaf1}) and (\ref{eqn:Mar-Vaf2}) by $R(\lambda; \tau;p)$
and $R(\lambda;\tau;p)^{\bullet}$ respectively.
Then the following two equivalent cut-and-join equations are satisfied:
\begin{eqnarray}
&& \frac{\partial R}{\partial \tau}
= \frac{\sqrt{-1}\lambda}{2} \sum_{i, j\geq 1} \left(ijp_{i+j}\frac{\partial^2R}{\partial p_i\partial p_j}
+ ijp_{i+j}\frac{\partial R}{\partial p_i}\frac{\partial R}{\partial p_j}
+ (i+j)p_ip_j\frac{\partial R}{\partial p_{i+j}}\right), \label{eqn:CutJoin1} \\
&& \frac{\partial R^{\bullet}}{\partial \tau}
= \frac{\sqrt{-1}\lambda}{2} \sum_{i, j\geq 1} \left(ijp_{i+j}\frac{\partial^2R^{\bullet}}{\partial p_i\partial p_j}
+ (i+j)p_ip_j\frac{\partial R^{\bullet}}{\partial p_{i+j}}\right). \label{eqn:CutJoin2}
\end{eqnarray}
\end{theorem}

\begin{proof}
We have by Proposition \ref{prop:CJ}
\begin{eqnarray*}
&& \frac{\partial R(\lambda;\tau;p)^{\bullet}}{\partial \tau} \\
& = & \sqrt{-1}\lambda \sum_{\mu, \nu}
\left( f_{\nu}(2) \frac{\chi_{\nu}(C(\mu))}{z_{\mu}} p_{\mu}\right)
e^{\sqrt{-1}(\tau+\frac{1}{2})\kappa_{\nu}\lambda/2}
V_{\nu}(\lambda)  \\
& = & \frac{\sqrt{-1}\lambda}{2}\sum_{i,j}\left(ijp_{i+j}\frac{\partial}{\partial p_i}\frac{\partial}{\partial p_j}
+ (i+j)p_ip_j \frac{\partial}{\partial p_{i+j}}\right) \sum_{\eta} \frac{\chi_{\nu}(\eta)}{z_{\eta}} p_{\eta}
e^{\sqrt{-1}(\tau+\frac{1}{2})\kappa_{\nu}\lambda/2}
V_{\nu}(\lambda) \\
& = & \frac{\sqrt{-1}\lambda}{2} \sum_{i, j\geq 1} \left(ijp_{i+j}
\frac{\partial^2R(\lambda;\tau;p)^{\bullet}}{\partial p_i\partial p_j}
+ (i+j)p_ip_j\frac{\partial R(\lambda;\tau;p)^{\bullet}}{\partial p_{i+j}}\right).
\end{eqnarray*}
(\ref{eqn:CutJoin1}) follows easily from (\ref{eqn:CutJoin2}).
\end{proof}

\subsection{The initial values}
It is possible to compute the initial values $\cC^{\bullet}(\lambda; 0; p)$ and
$R(\lambda; 0; p)$.
Indeed,
by Mumfords' relations:
$$\Lambda_g^{\vee}(1)\Lambda^{\vee}_g(-1) = (-1)^g,$$
one has
\begin{eqnarray*}
\cC(0, \lambda, p) & = & - \sum_{d > 0} \sqrt{-1}^{d+1}p_d \sum_{g \geq 0} \lambda^{2g-1}\int_{\Mbar_{g, 1}}
\frac{\Lambda_g^{\vee}(1)\Lambda_g^{\vee}(-1)\Lambda_g^{\vee}(0)}{1 - d\psi_1} \\
& = &  - \sum_{d > 0}  \sqrt{-1}^{d+1}p_d\sum_{g \geq 0} (d \lambda)^{2g-1}\int_{\Mbar_{g, 1}} \lambda_g\psi_1^{2g-2} \\
& = & - \sum_{d > 0} \frac{\sqrt{-1}^{d+1}p_d}{2d\sin (d\lambda/2)}.
\end{eqnarray*}
Here in the last equality we have used a result in \cite{Fab-Pan, Tia-Zho}.
Therefore,
\begin{eqnarray} \label{eqn:CInit}
&& \cC(\lambda; 0; p)^{\bullet} = \exp \left(-\sum_{d > 0} \frac{\sqrt{-1}^{d+1}p_d}{2d\sin (d\lambda/2)}\right).
\end{eqnarray}

\begin{theorem} \label{thm:RInit}
We have the following identity:
\begin{eqnarray}
&& \log\left(  \sum_{n \geq 0} \sum_{|\rho|=n}
\frac{e^{\frac{1}{4}\kappa_{\rho}\sqrt{-1}\lambda}}{\prod_{e \in \rho} 2\sin (h(e)\lambda/2)}
\frac{\chi_{\rho}(\eta)}{z_{\eta}} p_{\eta}\right)
= - \sum_{d \geq 1} \frac{\sqrt{-1}^{d+1}p_d}{2d\sin(d\lambda/2)}. \label{eqn:R}
\end{eqnarray}
\end{theorem}

To prove this result, we need some preliminary results.

\begin{lemma}
Introduce formal variables $x_1, \dots, x_n, \dots$ such that
$$p_i(x_1, \dots, x_n, \dots) = x_1^i + \cdots + x_n^i + \cdots.$$
Then for for any positive integer $n$,
we have
\begin{eqnarray} \label{eqn:pg}
&& \sum_{n \geq 0} t^n\sum_{|\rho|=n}
\frac{q^{n(\rho)}}{\prod_{e \in \rho} (1 - q^{h(e)})}
\frac{\chi_{\rho}(\eta)}{z_{\rho}} p_{\eta}
= \frac{1}{\prod_{i, j} (1 -tx_iq^{j-1})}.
\end{eqnarray}
\end{lemma}

\begin{proof}
Recall the following facts about Schur polynomials:
\begin{eqnarray}
&& s_{\rho}(x) = \sum_{\eta} \frac{\chi_{\rho}(\eta)}{z_{\eta}} p_{\eta}(x), \label{eqn:s} \\
&& s_{\rho}(1, q, q^2, \dots)
= \frac{q^{n(\rho)}}{\prod_{e \in \rho} (1 - q^{h(e)})}, \\
&& \sum_{n \geq 0} t^n \sum_{|\rho|=n} s_{\rho}(x)s_{\rho}(y)
= \frac{1}{\prod_{i, j} (1 -tx_iy_j)}.
\end{eqnarray}
Combining the last two identities,
one gets:
\begin{eqnarray*}
&& \sum_{n \geq 0} t^n \sum_{|\rho|=n}
\frac{q^{n(\rho)}}{\prod_{e \in \rho} (1 - q^{h(e)})}s_{\rho}(x)
= \frac{1}{\prod_{i,j} (1 - tx_iq^{j-1})}.
\end{eqnarray*}
The proof is completed by (\ref{eqn:s}).
\end{proof}

\begin{lemma}
For any partition $\rho$ we have
\begin{eqnarray}
&& \frac{1}{2} \sum_{e \in \rho} h(e) - n(\rho)
=  \frac{1}{4}\kappa_{\rho} + \frac{1}{2}|\rho|.
\end{eqnarray}
\end{lemma}

\begin{proof}
\begin{eqnarray*}
&& \frac{1}{2} \sum_{e \in \rho} h(e) - n(\rho) = \frac{1}{2} (n(\rho') - n(\rho) + |\rho|)  \\
& = & \frac{1}{2} (\sum_i \begin{pmatrix} \rho_i \\ 2 \end{pmatrix}
- \sum_i (i-1)\rho_i + |\rho|) \\
& = & \frac{1}{4} (\sum_i \rho_i(\rho_i - 1) - 2 \sum_i i\rho_i + 4|\rho|) \\
& = & \frac{1}{4}\kappa_{\rho} + \frac{1}{2}|\rho|.
\end{eqnarray*}
\end{proof}

We now come to the proof of Theorem \ref{thm:RInit}.
Let $q = e^{-\sqrt{-1}\lambda}$, and $t = \sqrt{-1}q^{1/2}$,
then we have
\begin{eqnarray*}
&& \sum_{n\geq 0} t^n\sum_{|\rho|=n}
\frac{q^{n(\rho)}}{\prod_{e \in \rho} (1 - q^{h(e)})}
\frac{\chi_{\rho}(\eta)}{z_{\rho}} p_{\eta} \\
& = & \sum_{n \geq 0} \sqrt{-1}^nq^{n/2} \sum_{|\rho|=n}
\frac{q^{n(\rho)-\frac{1}{2}\sum_{e\in \rho} h(e)}}{\prod_{e \in \rho} (q^{-h(e)/2} - q^{h(e)/2})}
\frac{\chi_{\rho}(\eta)}{z_{\rho}} p_{\eta} \\
& = & \sum_{n \geq 0} \sqrt{-1}^nq^{n/2} \sum_{|\rho|=n}
\frac{q^{-\frac{1}{4}\kappa_{\rho} - \frac{1}{2}n}}{\prod_{e \in \rho} (q^{-h(e)/2} - q^{h(e)/2})}
\frac{\chi_{\rho}(\eta)}{z_{\rho}} p_{\eta} \\
& = & \sum_{n \geq 0} \sum_{|\rho|=n}
\frac{e^{\frac{1}{4}\kappa_{\rho}\sqrt{-1}\lambda}}{\prod_{e \in \rho} 2\sin (h(e)\lambda/2)}
\frac{\chi_{\rho}(\eta)}{z_{\rho}} p_{\eta} \\
& = &  \sum_{n \geq 0} \sum_{|\rho|=n}
\frac{\chi_{\rho}(\eta)}{z_{\eta}} e^{\frac{1}{4}\kappa_{\rho}\lambda} V_{\rho} p_{\eta}.
\end{eqnarray*}
Hence by (\ref{eqn:pg}),
\begin{eqnarray*}
R(\lambda;0;p) & = & \log \left(  \sum_{n \geq 0} \sum_{|\rho|=n}
\frac{\chi_{\rho}(\eta)}{z_{\eta}} e^{\frac{1}{4}\kappa_{\rho}\lambda} V_{\rho} p_{\eta}\right)\\
& = & \log\left(  \sum_{n\geq 0} t^n\sum_{|\rho|=n}
\frac{q^{n(\rho)}}{\prod_{e \in \rho} (1 - q^{h(e)})}
\frac{\chi_{\rho}(\eta)}{z_{\rho}} p_{\eta}\right) \\
& = & \log  \frac{1}{\prod_{i, j} (1 -tx_iq^{j-1})}
= \sum_{i,j \geq 1} \sum_{d \geq 1} \frac{1}{d} t^d q^{d(j-1)} x_i^d \\
& = & \sum_{j \geq 1} \sum_{d \geq 1} \frac{1}{d} t^d q^{d(j-1)} p_d
= \sum_{d \geq 1} \frac{p_d}{d} \frac{t^d}{1 - q^d} \\
& = & - \sum_{d \geq 0} \frac{\sqrt{-1}^{d+1}p_d}{2d\sin(d\lambda/2)}.
\end{eqnarray*}

As a corollary
we prove the following result announced in \cite[Theorem 3]{Zho}.
It provides some evidence for Conjecture 2 and Conjecture 4 there.

\begin{theorem}
For any positive integer $n$ and any partition $\rho$ of $n$ we have
\begin{eqnarray} \label{eqn:Evidence}
&& \sum_{|\eta|=n} \frac{\chi_{\eta}(\rho)
e^{\sqrt{-1}\kappa_{\eta}\lambda/4}}
{\prod_{x \in \eta} 2\sin (h(x)\lambda/2)}
= \frac{\sqrt{-1}^{n - l(\rho)}}{2^{l(\rho)} \prod_k \sin^{m_k(\rho)}(k\lambda/2)}.
\end{eqnarray}
\end{theorem}

\begin{proof}
By (\ref{eqn:R}) we have
\begin{eqnarray*}
R(\lambda;0;p) & = & \sum_{n \geq 0} \sum_{|\rho|=n}
\frac{e^{\frac{1}{4}\kappa_{\rho}\sqrt{-1}\lambda}}{\prod_{e \in \rho} 2\sin (h(e)\lambda/2)}
\frac{\chi_{\rho}(\eta)}{z_{\eta}} p_{\eta} \\
& = & \exp \left(\sum_{k \geq 1} \frac{\sqrt{-1}^{k-1}p_k}{2k\sin(k\lambda/2)}\right)
= \prod_{k \geq 1} \sum_{m_k \geq 0} \left(\frac{\sqrt{-1}^{k-1}p_k}{2k\sin(k\lambda/2)}\right)^{m_k} \\
& = & \sum_{m_1, \dots, m_n \geq 0}
\frac{\sqrt{-1}^{\sum_k (k-1)m_k}}{2^{\sum_k m_k} \prod_k \sin^{m_k} (k\lambda/2)} \cdot
\prod_k \frac{p_k^{m_k}}{k^{m_k} m_k!} \\
& = & 1 + \sum_{n \geq 1} \sum_{|\rho|=n} \frac{\sqrt{-1}^{|\rho|-l(\rho)}}{2^{l(\rho)}
\prod_k \sin^{m_k(\rho)}(k\lambda/2)} \frac{p_{\rho}}{z_{\rho}}.
\end{eqnarray*}
The Theorem is proved by comparing the coefficients.
\end{proof}

\subsection{Low degree examples}

To illustrate the idea,
we consider some low degree examples.
By (\ref{eqn:CInit}) we have
\begin{eqnarray*}
\cC(\lambda; 0; p)^{\bullet} & = & \exp\left(
- \sum_{n > 0}  \frac{\sqrt{-1}^{n+1}p_n}{2n\sin (n\lambda/2)}\right) \\
& = & 1 + \frac{p_1}{2\sin (\lambda/2)}
+ \frac{\sqrt{-1}p_2}{4\sin \lambda} + \frac{p_1^2}{8\sin^2(\lambda/2)} + \cdots
\end{eqnarray*}

The degree $1$ case of (\ref{eqn:CutJoin4}) is
\begin{eqnarray*}
&& \frac{\partial \cC_{(1)}(\lambda; \tau)^{\bullet}}{\partial \lambda}
= 0.
\end{eqnarray*}
By the initial value:
$$\cC_{(1)}(\lambda; 0; p) = \frac{1}{2\sin (\lambda/2)},$$
we get
$$\cC_{(1)}(\lambda; \tau)^{\bullet} = \frac{1}{2\sin (\lambda/2)}.$$
This matches with the results in \cite{Zho}.
The degree $2$ case of (\ref{eqn:CutJoin4}) is
\begin{align*}
\frac{\partial \cC_{(2)}(\lambda;\tau)^{\bullet}}{\partial \tau}
& = \sqrt{-1}\lambda \cC_{(1^2)}(\lambda; \tau)^{\bullet}, &
\frac{\partial \cC_{(1^2)}(\lambda;\tau)^{\bullet}}{\partial \tau} & =
\sqrt{-1}\lambda \cC_{(2)}(\lambda; \tau)^{\bullet}.
\end{align*}
By the initial values
\begin{align*}
\cC_{(2)}(\lambda;0)^{\bullet} & =  \frac{\sqrt{-1}}{4\sin \lambda}, &
\cC_{(1^2)}(\lambda; 0)^{\bullet} & = \frac{1}{8\sin^2(\lambda/2)}
\end{align*}
one finds
\begin{align*}
\cC_{(2)}(\lambda;\tau)^{\bullet} &= \frac{\sqrt{-1}\cos (\tau\lambda)}{4\sin \lambda}
+ \frac{\sqrt{-1}\sin (\tau\lambda)}{8 \sin^2(\lambda/2)}
= \frac{\sqrt{-1}\sin[(\tau + \frac{1}{2})\lambda]}{4\sin(\lambda/2)\sin\lambda}, \\
\cC_{(1^2)}(\lambda;\tau)^{\bullet} & = - \frac{\sin (\tau\lambda)}{4\sin \lambda}
+ \frac{\cos (\tau\lambda)}{8\sin^2(\lambda/2)}
=  \frac{\cos[(\tau + \frac{1}{2})\lambda]}{4\sin(\lambda/2)\sin\lambda}.
\end{align*}
Hence we have
\begin{eqnarray*}
&& \cC(\lambda; \tau, p) \\
& = & \log \left( 1 + \frac{p_1}{2\sin (\lambda/2)}
+ p_2\frac{\sqrt{-1}\sin[(\tau + \frac{1}{2})\lambda]}{4\sin(\lambda/2)\sin\lambda}
+ p_1^2 \frac{\cos[(\tau + \frac{1}{2})\lambda]}{4\sin(\lambda/2)\sin\lambda} + \cdots \right) \\
& = & 1 + \frac{p_1}{2\sin (\lambda/2)}
+ p_2\frac{\sqrt{-1}\sin[(\tau + \frac{1}{2})\lambda]}{4\sin(\lambda/2)\sin\lambda}
+ p_1^2 \left(\frac{\cos[(\tau + \frac{1}{2})\lambda]}{4\sin(\lambda/2)\sin\lambda}
- \frac{1}{8\sin^2 (\lambda/2)}\right) \\
&& + \cdots  \\
& = & 1 + \frac{p_1}{2\sin (\lambda/2)}
+ p_2\frac{\sqrt{-1}\sin[(\tau + \frac{1}{2})\lambda]}{4\sin(\lambda/2)\sin\lambda}
- p_1^2 \frac{\sin(\tau\lambda/2)\sin[(\tau + 1)\lambda/2]}{4\sin(\lambda/2)\sin\lambda} + \cdots.
\end{eqnarray*}
This matches with the results proved in \cite{Zho}.

\subsection{Some consequence of the cut-and-join equation}
The following result is announced in \cite[Theorem 4]{Zho}.
It leads to a simple proof of the $\lambda_g$ conjecture.
See \cite{LLZ3}.

\begin{theorem}
One has for partition $\mu$ of $d$,
\begin{equation}
\begin{split}
& \lim_{\tau \to 0}  \lambda^{2-l(\mu)} \cdot\frac{1}{(\tau(\tau+1))^{l(\mu)-1}}  \cdot
\prod_{i=1}^{l(\mu)} \frac{\mu_i \cdot \mu_i!}{\prod_{j=1}^{\mu_i-1}(j+\mu_i\tau)} \cdot
\frac{\prod_j m_j(\mu)!}{\sqrt{-1}^{|\mu|+l(\mu)}}\\
& \cdot \sum_{n \geq 1} \frac{(-1)^n}{n} \sum_{\cup_{i=1}^n \mu^i = \mu}
\sum_{|\nu^i|=|\mu^i|} \prod_{i=1}^n
\frac{\chi_{\nu^i}(C(\mu^i))}{z_{\mu^i}}
\cdot e^{\sqrt{-1}(\tau+\frac{1}{2})\kappa_{\nu^i}\lambda/2}
\cdot V_{\nu^i}(\lambda) \\
= & \prod_{i=1}^{l(\mu)} \mu_i^2 \cdot d^{l(\mu)-3} \cdot \frac{d\lambda/2}{\sin (d\lambda/2)}.
\end{split} \end{equation}
\end{theorem}

\begin{proof}
Write
$$R(\lambda;\tau;p) = \sum_{\mu} R_{\mu}(\lambda;\tau)p_{\mu}.$$
By (\ref{eqn:R}) we have
\begin{eqnarray}
&& R_{\mu}(\lambda;0) = \delta_{l(\mu), 1} \frac{\sqrt{-1}^{|\mu| - 1}}{2|\mu|\sin(|\mu|\lambda/2)}.
\label{eqn:R0}
\end{eqnarray}
Since
\begin{eqnarray*}
&& \lim_{\tau \to 0}  \lambda^{2-l(\mu)} \cdot\frac{1}{(\tau+1)^{l(\mu)-1}}  \cdot
\prod_{i=1}^{l(\mu)} \frac{\mu_i \cdot \mu_i!}{\prod_{j=1}^{\mu_i-1}(j+\mu_i\tau)} \cdot
\frac{\prod_j m_j(\mu)!}{\sqrt{-1}^{|\mu|+l(\mu)}} \\
& = &  \lambda^{2-l(\mu)}  \cdot
\prod_{i=1}^{l(\mu)} \mu_i^2  \cdot
\frac{\prod_j m_j(\mu)!}{\sqrt{-1}^{|\mu|+l(\mu)}},
\end{eqnarray*}
it suffices to prove
\begin{eqnarray} \label{eqn:Limit}
&& \lim_{\tau \to 0} \frac{R_{\mu}(\lambda;\tau)}{\tau^{l(\mu)-1}}
= -\frac{\sqrt{-1}^{|\mu|+l(\mu)}}{\prod_j m_j(\mu)!}\lambda^{l(\mu)-1}
\frac{|\mu|^{l(\mu)-2}}{2\sin (2|\mu|\lambda/2)}.
\end{eqnarray}
When $l(\mu) = 1$,
(\ref{eqn:Limit}) is just (\ref{eqn:R0}).
When $l(\mu) > 1$ we will prove (\ref{eqn:Limit}) by L'Hospital's rule.
We will use the following notations:
\begin{align*}
CF & = \frac{1}{2} \sum_{i,j} (i+j)p_ip_j \frac{\partial F}{\partial p_{i+j}}, \\
JF & = \frac{1}{2} \sum_{i,j} ijp_{i+j}\frac{\partial^2 F}{\partial p_i\partial p_j}.
\end{align*}
Then by (\ref{eqn:CutJoin2}) and (\ref{eqn:R0}),
\begin{eqnarray*}
&& \left. \frac{\partial^k}{\partial \tau^k} R(\lambda;\tau;p)\right|_{\tau = 0}
= (\sqrt{-1}\lambda)^k C^kR(\lambda;0;p) + \cdots  \\
& = & (\sqrt{-1}\lambda)^k C^k\sum_{d \geq 1} \frac{\sqrt{-1}^{d-1}p_d}{2d\sin(d\lambda/2)}
+ \cdots,
\end{eqnarray*}
where $\cdots$ stands for terms that involve the join operator $J$ at least once.
Note the right-hand side can be written as a linear combination of $p_{\nu}$ with $l(\nu) \leq k+1$.
Furthermore those terms with $l(\nu) = k+1$ are all obtained from cuts for $k$ times.
Hence by comparing the coefficients of $p_{\mu}$ for $l(\mu) = l$ on both sides,
one sees that for $k < l-1$,
one has
$$\left. \frac{\partial^k}{\partial \tau^k} R_{\mu}(\lambda;\tau)\right|_{\tau = 0} = 0,$$
and by Lemma \ref{lm:CP} below,
\begin{eqnarray*}
&& \sum_{l(\mu) = l} \left. \frac{\partial^{l-1}}{\partial \tau^{l-1}} R_{\mu}(\lambda;\tau)\right|_{\tau = 0}p_{\mu}
= (\sqrt{-1}\lambda)^{l-1} C^{l-1}\sum_{d \geq 1} \frac{\sqrt{-1}^{d-1}p_d}{2d\sin(d\lambda/2)} \\
& = & \lambda^{l-1}
\sum_{d \geq 1} \sqrt{-1}^{l+d-2} \frac{C^{l-1}p_d}{2d\sin(d\lambda/2)} \\
& = & \lambda^{l-1}
\sum_{d \geq 1} \sqrt{-1}^{l+d-2} \frac{1}{2d\sin(d\lambda/2)}
\sum_{|\mu|=d, l(\mu) =l} \frac{(l-1)!d^{l-1}}{\prod_{j} m_j(\mu)!}  p_{\mu} \\
& = & \sum_{l(\mu)=l} \lambda^{l-1}
 \sqrt{-1}^{l+|\mu|-2} \frac{1}{2|\mu|\sin(|\mu|\lambda/2)}
\frac{(l-1)!|\mu|^{l-1}}{\prod_{j} m_j(\mu)!}  p_{\mu}.
\end{eqnarray*}
Therefore,
(\ref{eqn:Limit}) follows by L'Hospital's rule.
\end{proof}

\begin{lemma} \label{lm:CP}
For $l \geq 1$ we have
\begin{eqnarray} \label{eqn:CP}
&& C^{l-1} p_d = \sum_{|\mu|=d, l(\mu) = l} \frac{(l-1)!d^{l-1}}{\prod_{j} m_j(\mu)!}  p_{\mu}
\end{eqnarray}
\end{lemma}

\begin{proof}
We prove (\ref{eqn:CP}) by induction on $l$.
When $l=1$,
it is trivial.
Suppose it holds for $l=k$.
Then we have
\begin{eqnarray*}
&& C^{k} p_d = \sum_{|\mu|=d, l(\mu) = k} \frac{(k-1)!d^{k-1}}{\prod_{j} m_j(\mu)!}  Cp_{\mu} \\
&= & \sum_{|\mu|=d, l(\mu) = k} \frac{(k-1)!d^{k-1}}{\prod_{j} m_j(\mu)!}
\frac{1}{2} \sum_{i,j}(i+j) m_{i+j}(\mu) p_{\nu} \;\;\;\;\;
(\nu \in C_{i,j}(\mu)) p_{\nu} \\
& = & \sum_{|\nu|=d, l(\nu) = k+1} \frac{(k-1)!d^{k-1}}{2\prod_j m_j(\nu)!}
\left(\sum_{i\neq j} m_i(\nu)m_j(\nu) (i+j) + \sum_i m_i(\nu)(m_i(\nu)-1) \cdot 2i\right)\\
& = & \sum_{|\nu|=d, l(\nu) = k+1} \frac{(k-1)!d^{k-1}}{\prod_j m_j(\nu)!}
\left(\sum_{i, j} m_i(\nu)m_j(\nu) (i+j) - \sum_i m_i(\nu) \cdot i\right) p_{\nu}\\
& = & \sum_{|\nu|=d, l(\nu) = k+1} \frac{(k-1)!d^{k-1}}{\prod_j m_j(\nu)!}
\sum_i  m_i(\nu)i \left( \sum_j m_j(\nu) -1\right) p_{\nu} \\
& = &  \sum_{|\nu|=d, l(\nu) = k+1} \frac{k!d^{k}}{\prod_j m_j(\nu)!} p_{\nu}.
\end{eqnarray*}
This finishes the proof.
\end{proof}

\section{Related Results for Hurwitz Numbers}
\label{sec:Hurwitz}

In this section we present some results for Hurwitz numbers analogous to the above results for
conifold Hodge integrals.
We also get some closed formulas for some generating series for Hurwitz numbers and related Hodge integrals
by the cut-and-join equations.

\subsection{Hurwitz numbers and Hodge integrals}

Let $X$ be a Riemann surface of genus $h$.
Given $n$ partitions $\eta^1, \dots, \eta^n$ of $d$,
denote by $H^X_d(\eta^1, \dots, \eta^n)^{\bullet}$ and $H^X_d(\eta^1, \dots, \eta^n)^{\circ}$
the weight counts of possibly disconnected and connected Hurwitz covers of type $(\eta^1, \dots, \eta^n)$
respectively.
We will use the following formula for Hurwitz numbers (see e.g. \cite{Dij}):
\begin{eqnarray} \label{eqn:Burnside}
&& H^X_d(\eta^1, \dots, \eta^n)^{\bullet}
= \sum_{\rho \vdash d} \left(\frac{\dim R_\rho}{d!}\right)^{2 - 2h}
\prod_{i=1}^n |C_{\eta^i}|\frac{\chi_{\rho}(\eta^i)}{\dim R_{\rho}}.
\end{eqnarray}
It is sometimes referred to the Burnside formula.

Given a partition $\eta$ of length $l$ and genus $g$,
let
$$r = 2g -2 + d + l - 2d h.$$
Then by Riemann-Hurwitz formula,
an almost simple Hurwitz cover of degree $d$ with ramification type $\eta$ at the possibly nonsimple
ramification point and $r$ other simple ramification points should have genus $g$.
Given a partition $\eta$,
denote
\begin{eqnarray*}
&& H_h^{g, d}({\eta})^{\circ} = H^X_d(\eta^1, \dots, \eta^{r+1})^{\circ}, \\
&& H_h^{g, d}({\eta})^{\bullet} = H^X_d(\eta^1, \dots, \eta^{r+1})^{\bullet},
\end{eqnarray*}
for $\eta^1 = \cdots = \eta^r = (2)$.
We have by (\ref{eqn:Burnside}):
\begin{eqnarray} \label{eqn:Simple}
H^g_h(\eta)^{\bullet}
& = & \sum_{|\rho| = d}  \left(\frac{\dim R_\rho}{d!}\right)^{2-2h}
\cdot f_\rho(2)^{2g-2+d+l(\eta)-2|\eta|h} \cdot f_\rho(\eta).
\end{eqnarray}

The ELSV formula \cite{ELSV, Gra-Vak} provides a deep connection between Hurwitz numbers and Hodge integrals:
\begin{eqnarray} \label{eqn:ELSV}
H_0^g(\eta)^{\circ} = \frac{(2g-2+|\eta|+l(\eta))!}{|\Aut (\eta)|}
\prod_{i=1}^{l(\eta)} \frac{\eta_i^{\eta_i}}{\eta_i!}
\int_{\Mbar_{g, l(\eta)}} \frac{\Lambda_g^{\vee}(1)}{\prod_{i=1}^{l(\eta)} (1 - \eta_i \psi_i)},
\end{eqnarray}
where
$$|\Aut (\eta)| = \prod_j m_j(\eta)!.$$
One can use this formula to transfer the results on Hurwitz numbers below to results on Hodge integrals.

\subsection{Generating functions of almost simple Hurwitz numbers}
Consider
\begin{eqnarray*}
&& \Phi_h^{\eta}(\lambda)^{\circ}
= \sum_{g \geq 0}
H^g_h(\eta)^{\circ}  \frac{\lambda^{2g-2+|\eta|+l(\eta)-2|\eta|h}}{(2g-2 + |\eta| + l(\eta)-2|\eta|h)!}, \\
&& \Phi_h^{\eta}(\lambda)^{\bullet}
= \sum_{g \geq 0} H^g_h(\eta)^{\bullet}
\frac{\lambda^{2g-2+|\eta|+l(\eta)-2|\eta|h}}{(2g-2 + |\eta| + l(\eta)-2|\eta|h)!}, \\
&& \Phi_h(\lambda, p)^{\circ} = \sum_{\eta} \Phi_h^{\eta}(\lambda)^{\circ} \cdot p_{\eta}, \\
&& \Phi_h(\lambda, p)^{\bullet}= 1+\sum_{\eta} \Phi^{\eta}_h(\lambda)^{\bullet} \cdot p_{\eta}.
\end{eqnarray*}
The usual relationship between connected and disconnected Hurwitz numbers
is (see e.g. \cite{Gou-Jac}):
\begin{eqnarray} \label{eqn:Conn}
&&\Phi_h(\lambda, p)^{\circ} = \log \Phi_h(\lambda, p)^{\bullet}.
\end{eqnarray}

\begin{theorem} \label{thm:Hurwitz}
We have
\begin{eqnarray}
\Phi_h(\lambda,p)^{\bullet}
& = &1+ \sum_{d \geq 1} \sum_{|\eta|= d} U_h^{\eta}(\lambda) \cdot p_{\eta}, \label{eqn:PhiBullet}\\
\Phi_h(\lambda,p)^{\circ}
& = & \sum_{ n \geq 1} \frac{(-1)^{n-1}}{n} \sum_{\eta}
\sum_{\cup_{j=1}^n \eta_j = \eta} \prod_{j=1}^n U_h^{\eta_j} \cdot p_{\eta}, \label{eqn:PhiCirc}
\end{eqnarray}
where
\begin{eqnarray} \label{eqn:W}
&& U_h^{\eta}(\lambda) =  \sum_{|\rho|=|\eta|} \left(\frac{\dim R_{\rho}}{|\eta|!}\right)^{2-2h} \cdot
\exp \left[f_{\rho}(2)\lambda\right] \cdot f_{\rho}(\eta).
\end{eqnarray}
\end{theorem}

\begin{proof}
The first identity is an easy consequence of (\ref{eqn:Burnside}).
The second identity is from the first identity by taking logarithm.
\end{proof}

\begin{remark}
By Proposition \ref{prop:f},
one has
\begin{eqnarray*}
&& U_h^{\eta}(\lambda)\\
& = & \frac{1}{2} \sum_{|\rho|=|\eta|} \left[\left(\frac{\dim R_{\rho}}{|\eta|!}\right)^{2-2h} \cdot
e^{f_{\rho}(2)\lambda}\cdot f_{\rho}(\eta)
+  \left(\frac{\dim R_{\rho^t}}{|\eta|!}\right)^{2-2h} \cdot
e^{f_{\rho^t}(2)\lambda} \cdot f_{\rho^t}(\eta)\right] \\
& = &   \sum_{|\rho|=|\eta|} \left(\frac{\dim R_{\rho}}{|\eta|!}\right)^{2-2h} \cdot
\frac{e^{f_{\rho}(2)\lambda}
+ (-1)^{|\eta|-l(\eta)}e^{f_{\rho}(2)\lambda}}{2} \cdot f_{\rho}(\eta).
\end{eqnarray*}
Hence depending on the parity of $|\eta|-l(\eta)$,
$U_h^{\eta}(\lambda)$ is a linear combination of $\cosh$ or $\sinh$ functions.
\end{remark}

\subsection{Generating functions of simple Hurwitz numbers}

We now consider the case of simple Hurwitz numbers.
Introduce
\begin{eqnarray*}
&& \Phi_h(\lambda, q)^{\circ}_s
= \sum_{d>0, g \geq 0}
H^g_h(1^d)^{\circ} \cdot \frac{\lambda^{2g-2+2d-2dh}}{(2g-2 + 2d-2dh)!} \cdot q^d, \\
&& \Phi_h(\lambda, q)^{\bullet}_s
= 1+\sum_{d>0, g \geq 0}
H^g_h(1^d)^{\bullet} \cdot \frac{\lambda^{2g-2+2d-2dh}}{(2g-2 + 2d-2dh)!} \cdot q^d.
 \end{eqnarray*}
They are obtained from $\Phi_h(\lambda, p)^{\circ}$ and $\Phi_h(\lambda, p)^{\bullet}$, respectively,
by taking $p_1 =q$, $p_2 = \cdots =p_n = \cdots = 0$.
By noting
\begin{eqnarray*}
&& f_{\rho}(1^d) = 1, \\
&& f_{\rho}(1^{d-2}2) = \begin{pmatrix} d \\2\end{pmatrix} \frac{\chi_{\rho}(2)}{\dim R_{\rho}}.
\end{eqnarray*}
one easily gets
$$\Phi(\lambda, q)^{\circ}_s = \log \Phi(\lambda, q)^{\bullet}_s,$$
and

\begin{theorem}
We have
\begin{eqnarray}
\Phi_h(\lambda,q)^{\bullet}_s
& = & 1+ \sum_{d \geq 1} \sum_{|\rho| = d}  \left(\frac{\dim R_\rho}{d!}\right)^{2-2h} \cdot
\exp \left[\begin{pmatrix} d \\2\end{pmatrix}
\frac{\chi_{\rho}(2)}{\dim R_{\rho}}\lambda\right] \cdot q^d  \label{eqn:Burnside1} \\
\Phi_h(\lambda,q)^{\circ}_s & = & \sum_{ n \geq 1} \frac{(-1)^{n-1}}{n} \sum_{d > 0} q^d
\sum_{\sum_{j=1}^n |\rho_j| = d}
\prod_{j=1}^n \left( \frac{\dim R_{\rho_j}}{|\rho_j|!}\right)^{2-2h} \label{eqn:Burnside2} \\
&& \cdot
\exp \left[\begin{pmatrix} |\rho_j| \\2\end{pmatrix}
\frac{\chi_{\rho_j}(2)}{\dim R_{\rho_j}}\lambda\right]. \nonumber
\end{eqnarray}
\end{theorem}

\begin{remark}
Since $\Phi_h(\lambda, q)^{\bullet}$ is an even function of $\lambda$,
one easily gets
\begin{eqnarray*}
\Phi_h(\lambda,q)^{\bullet}_s
& = & 1+ \sum_{d \geq 1} \sum_{|\rho| = d}  \left(\frac{\dim R_\rho}{d!}\right)^{2-2h} \cdot
\cosh \left[\begin{pmatrix} d \\2\end{pmatrix}
\frac{\chi_{\rho}(2)}{\dim R_{\rho}}\lambda\right] \cdot q^d.
\end{eqnarray*}
This is proved in \cite{Mon-Son-Son} by a different method.
\end{remark}

\subsection{Low degree examples in genus $0$}
\label{sec:Example01}

Using the character tables
one easily gets
\begin{eqnarray*}
&& \Phi_0(\lambda, p)^{\bullet} \\
& = & 1 + p_1 + \frac{1}{4}(e^t - e^{-t})p_2 + \frac{1}{4} (e^t + e^{-t})p_1^2 \\
& + & \frac{1}{36}(2e^{3\lambda}-4+2e^{-3\lambda})p_3
+ \frac{1}{36}(3e^{3\lambda}-3e^{-3\lambda})p_1p_2
+ \frac{1}{36}(e^{3\lambda}+4+e^{-3\lambda})p_1^3 \\
& + & \frac{1}{576}(6e^{6\lambda}+ 18e^{2\lambda} + 18e^{-2\lambda}-6e^{-6\lambda})p_4
+  \frac{1}{576}(8e^{6\lambda}-16+8e^{-6\lambda})p_1p_3 \\
& + & \frac{1}{576}(3e^{6\lambda}-9e^{2\lambda} +12 + 9e^{-2\lambda}+3e^{-6\lambda})p_2^2 \\
&+ &  \frac{1}{576}(6e^{6\lambda}+18e^{2\lambda} - 18e^{-2\lambda}-6e^{-6\lambda})p_1^2p_2\\
&+& \frac{1}{576}(e^{6\lambda}+9e^{2\lambda} + 4 + 9e^{-2\lambda}+e^{-6\lambda})p_1^4
+ \cdots \\
& = & 1 + p_1 + \frac{1}{2}p_2 \sinh \lambda + \frac{1}{2} p_1^2 \cosh \lambda
+ \frac{1}{9}p_3(\cosh (3\lambda) -1) \\
&+& \frac{1}{6}p_1p_2\sinh (3\lambda)
+ \frac{1}{18}p_1^3(\cosh (3\lambda) + 2) \\
&+&  \frac{1}{48}(\sinh (6\lambda)  - 3 \sin (2\lambda)) p_4
+\frac{1}{36}(\cosh (6\lambda)-1)p_1p_3 \\
&+& \frac{1}{96}(\cosh (6\lambda) - 3\cosh (2\lambda) + 2)p_2^2
 + \frac{1}{48}(\sinh (6\lambda) + 3 \sinh (2\lambda))p_1^2p_2 \\
&+&  \frac{1}{288}(\cosh(6\lambda) + 9\cosh(2\lambda) + 2) p_1^4+ \cdots
\end{eqnarray*}
After taking logarithm and some simple algebraic manipulations,
one gets
\begin{eqnarray*}
&& \Phi_0(\lambda, p)^{\circ} \\
& = & p_1 + \frac{1}{2}p_2\sinh \lambda
+  p_1^2\sinh^2(\lambda/2) + \frac{2}{9}p_3\sinh^2(3\lambda/2)\\
& + & p_1p_2\left(\frac{1}{6}\sinh(3\lambda) - \frac{1}{2}\sinh \lambda\right)
+ p_1^3\left(\frac{1}{9}\sinh^2(3\lambda/2) - \sinh^2(\lambda/2)\right) + \cdots
\end{eqnarray*}

\subsection{Low degree examples in genus $1$}
\label{sec:Example11}

Similarly one has
\begin{eqnarray*}
&& \Phi_1(\lambda, p)^{\bullet} \\
& = & 1 + p_1 + (e^t - e^{-t})p_2 + \frac{1}{4} (e^t + e^{-t})p_1^2
+  (e^{3\lambda}-1+e^{-3\lambda})p_3 \\
& + & (3e^{3\lambda}-3e^{-3\lambda})p_1p_2
+ (e^{3\lambda}+1+e^{-3\lambda})p_1^3 + \cdots \\
& = & 1 + p_1 + 2p_2 \sinh \lambda + 2 p_1^2 \cosh \lambda
+ p_3(4\cosh (3\lambda) -1) + 6p_1p_2\sinh (3\lambda) \\
& + & p_1^3(2\cosh (3\lambda) + 1)+
\cdots
\end{eqnarray*}
Taking logarithm,
one gets
\begin{eqnarray*}
&& \Phi_1(\lambda, p)^{\circ} \\
& = & p_1 + 2p_2\sinh \lambda
+  p_1^2(2\cosh\lambda - \frac{1}{2}) + p_3(4\sinh(3\lambda)-1)\\
& + & p_1p_2\left(6\sinh(3\lambda) - 2\sinh \lambda\right)
+ p_1^3\left(2\cosh(3\lambda) - 2\cosh\lambda+\frac{4}{3}\right) + \cdots
\end{eqnarray*}

\subsection{The cut-and-join equations for Hurwtiz numbers}

By combining Theorem \ref{thm:Hurwitz} and Proposition \ref{prop:CJ},
one easily obtains the following:

\begin{proposition} \label{prop:CutJoin}
We have the following equations:
\begin{eqnarray}
&& \frac{\partial \Phi_h^{\bullet}}{\partial \lambda}
= \frac{1}{2} \sum_{i, j\geq 1} \left(ijp_{i+j}\frac{\partial^2\Phi_h^{\bullet}}{\partial p_i\partial p_j}
+ (i+j)p_ip_j\frac{\partial \Phi_h^{\bullet}}{\partial p_{i+j}}\right), \label{eqn:CutJoin3} \\
&& \frac{\partial \Phi_h^{\circ}}{\partial \lambda}
= \frac{1}{2} \sum_{i, j\geq 1} \left(ijp_{i+j}\frac{\partial^2\Phi_h^{\circ}}{\partial p_i\partial p_j}
+ ijp_{i+j}\frac{\partial \Phi_h^{\circ}}{\partial p_i}\frac{\partial \Phi_h^{\circ}}{\partial p_j}
+ (i+j)p_ip_j\frac{\partial \Phi_h^{\circ}}{\partial p_{i+j}}\right).  \label{eqn:CutJoin4}
\end{eqnarray}
\end{proposition}

\begin{remark}
Equation (\ref{eqn:CutJoin4}) was proved in \cite{Gou-Jac-Vai} by combinatorial method for $h=0$,
and in \cite{Li-Zha-Zhe} by symplectic method for general $h$,
where the original forms of the equations has extra variables $z$ and $x$.
The present form of the equation is obtained by taking $z=x=1$.
Equation (\ref{eqn:CutJoin3}) can be obtained from (\ref{eqn:CutJoin4}) by (\ref{eqn:Conn}).
\end{remark}

Again, by writing (\ref{eqn:CutJoin3}) as a sequence of systems of ODEs,
one sees that
$\{\Phi_h^{\eta}(\lambda)^{\bullet}: |\eta|=d\}$ is determined by
the initial values $\{\Phi_h^{\eta}(0)^{\bullet}: |\eta|=d\}$.

In the above we have used (\ref{eqn:PhiBullet}) and (\ref{eqn:PhiCirc}) to prove the Hurwitz numbers satisfy
the cut-and-join equations.
We can actually reverse the procedure to show Hurwitz numbers satisfy (\ref{eqn:PhiBullet}) and (\ref{eqn:PhiCirc})
as follows.
This inspires the approach to the proof of Mari\~{n}o-Vafa formula described in last section.
The right-hand side of (\ref{eqn:PhiBullet}) satisfies (\ref{eqn:CutJoin3})
by  our combinatorial result Proposition \ref{prop:CJ};
the left-hand side satisfies the same equation by the geometric approach of \cite{Li-Zha-Zhe}.
Hence it suffices to check that both sides have the same initial values.
On the one hand,
it is easy to see that
\begin{eqnarray}
&& \Phi^{\eta}_0(0)^{\bullet} = \begin{cases}
0, & \eta \neq (1^{|\eta|}), \\
\frac{1}{|\eta|!}, & \eta = (1^{|\eta|}),
\end{cases} \label{eqn:Init0}
\end{eqnarray}
i.e.,
$$\Phi_0(0, p)^{\bullet} = e^{p_1}.$$
On the other hand,
by taking $\lambda = 0$ on the right-hand side of (\ref{eqn:PhiBullet}),
one gets:
\begin{eqnarray*}
&& \sum_{|\rho|=|\eta|}  \left(\frac{|\eta|!}{\dim R_{\rho}}\right)^{-2} f_{\rho}(\eta)
= \sum_{|\rho|=|\eta|}  |C_{\eta}|\frac{\dim R_{\rho} \cdot \chi_{\rho}(\eta)}{(|\eta|!)^2} \\
& = & \frac{|C_{\eta|}}{|\eta|!} \chi_{L}(g)
= \frac{|C_{\eta}|}{(|\eta|!)^2} |\eta|!\delta_{\eta, (1^{|\eta|})}
= \frac{1}{|\eta|!}\delta_{\eta, (1^{|\eta|})},
\end{eqnarray*}
where $L$ is the left regular representation.
This matches with (\ref{eqn:Init0}).
Hence by Proposition \ref{prop:CutJoin} $\Phi_0(\lambda, p)^{\bullet}$ is given by (\ref{eqn:PhiBullet}).
In genus $h$,
(\ref{eqn:PhiBullet}) predicts
\begin{eqnarray} \label{eqn:Initial}
\Phi_h^{\eta}(0)^{\bullet}
& = &  U_h^{\eta}(0)
 =  \sum_{|\rho|=|\eta|} \left(\frac{\dim R_{\rho}}{|\eta|!}\right)^{2-2h} \cdot f_{\rho}(\eta).
\end{eqnarray}
If one can establish this directly,
then one recovers (\ref{eqn:PhiBullet}) by Proposition \ref{prop:CutJoin}.

By the method of \cite{LLZ2} one can show that the Hodge integrals on
the right-hand side of
the ELSV formula satisfies the cut-and-join equation,
and it is also easy to see that both sides have the same initial values.
This leads to a proof of the ELSV formula in the same fashion of
the proof of the Mari\~no-Vafa formula presented in \cite{LLZ1, LLZ2}.
See \cite{LLZ3} for details.

\subsection{Low degree examples in genus $0$}

The degree $1$ case of (\ref{eqn:CutJoin3}) is
\begin{eqnarray*}
&& \frac{\partial \Phi_h^{(1)}(\lambda)^{\bullet}}{\partial \lambda} = 0.
\end{eqnarray*}
This is compatible with the fact that
$$\mu^g_h(1)^{\bullet} = \delta_{g,h}.$$

The degree $2$ case of (\ref{eqn:CutJoin3}) is
\begin{align*}
\frac{\partial \Phi^{(2)}_h(\lambda)^{\bullet}}{\partial \lambda} &= \Phi^{(11)}_h(\lambda)^{\bullet},
& \frac{\partial \Phi^{(1^2)}_h(\lambda)^{\bullet}}{\partial \lambda} & = \Phi^{(2)}_h(\lambda)^{\bullet}.
\end{align*}
By the initial values
\begin{align*}
\Phi^{(2)}_0(0)^{\bullet} &=0,
& \Phi^{(1^2)}_0(0)^{\bullet} & = \frac{1}{2},
\end{align*}
one finds
\begin{align*}
\Phi^{(2)}_0(\lambda)^{\bullet} &= \frac{1}{2}\sinh \lambda,
& \Phi^{(1^2)}_0(\lambda)^{\bullet} & = \frac{1}{2}\cosh \lambda.
\end{align*}

The degree $3$ case of (\ref{eqn:CutJoin3}) is
\begin{align*}
\frac{\partial \Phi^{(3)}_h(\lambda)^{\bullet}}{\partial \lambda} & = 2\Phi^{(12)}_h(\lambda)^{\bullet}, \\
\frac{\partial \Phi^{(12)}_h(\lambda)^{\bullet}}{\partial \lambda} & = 3 \Phi^{(3)}_h(\lambda)^{\bullet}
+ 3 \Phi^{(1^3)}(\lambda)^{\bullet}, \\
\frac{\partial \Phi^{(1^3)}_h(\lambda)^{\bullet}}{\partial \lambda} & = \Phi^{(12)}_h(\lambda)^{\bullet}.
\end{align*}
By the initial values
\begin{align*}
\Phi^{(3)}_0(0)^{\bullet} &=0, &
\Phi^{(12)}_0(0)^{\bullet} &=0,
& \Phi^{(1^3)}_0(0)^{\bullet} & = \frac{1}{6},
\end{align*}
one finds
\begin{align*}
\Phi^{(3)}_0(\lambda)^{\bullet} &=  \frac{1}{9}(\cosh (3\lambda) -1) ,
& \Phi^{(12)}_0(\lambda)^{\bullet} & = \frac{1}{6}\sinh (3\lambda),
& \Phi^{(1^3)}_0(\lambda)^{\bullet} & = \frac{1}{18}(\cosh (3\lambda) + 2).
\end{align*}

The degree $4$ case of (\ref{eqn:CutJoin3}) is
\begin{align*}
\frac{\partial \Phi^{(4)}_h(\lambda)^{\bullet}}{\partial \lambda}
& = 3\Phi^{(13)}_h(\lambda)^{\bullet} + 4 \Phi^{(2^2)}_h(\lambda)^{\bullet}, \\
\frac{\partial \Phi^{(13)}_h(\lambda)^{\bullet}}{\partial \lambda}
& = 4 \Phi^{(4)}_h(\lambda)^{\bullet} + 4 \Phi^{(1^22)}(\lambda)^{\bullet}, \\
\frac{\partial \Phi^{(2^2)}_h(\lambda)^{\bullet}}{\partial \lambda}
& = 2\Phi^{(4)}_h(\lambda)^{\bullet} + \Phi^{(1^22)}_h(\lambda)^{\bullet}, \\
\frac{\partial \Phi^{(1^22)}_h(\lambda)^{\bullet}}{\partial \lambda}
& = 2 \Phi^{(2^2)}_h(\lambda)^{\bullet} + 3\Phi^{(13)}_h(\lambda)^{\bullet}+6\Phi^{(1^4)}_h(\lambda)^{\bullet}, \\
\frac{\partial \Phi^{(1^4)}_h(\lambda)^{\bullet}}{\partial \lambda}
& = \Phi^{(1^22)}_h(\lambda)^{\bullet}.
\end{align*}
By the initial values (\ref{eqn:Init0})
one finds
\begin{align*}
\Phi^{(4)}_0(\lambda)^{\bullet} &=  \frac{1}{48}(\sinh (6\lambda)  - 3 \sin (2\lambda)) ,\\
\Phi^{(13)}_0(\lambda)^{\bullet} & = \frac{1}{36}(\cosh (6\lambda)-1), \\
\Phi^{(2^2)}_0(\lambda)^{\bullet} & = \frac{1}{96}(\cosh (6\lambda) - 3\cosh (2\lambda) + 2), \\
\Phi^{(1^22)}_0(\lambda)^{\bullet} & =  \frac{1}{48}(\sinh (6\lambda) + 3 \sinh (2\lambda)), \\
\Phi^{(1^4)}_0(\lambda)^{\bullet} & = \frac{1}{288}(\cosh(6\lambda) + 9\cosh(2\lambda) + 2).
\end{align*}
These results match with calculations in \S \ref{sec:Example01}.

\subsection{Low degree examples in genus $1$}
Note $\Phi_1^{\eta}(\lambda)^{\bullet}$ satisfy the same equations as $\Phi_0^{\eta}(\lambda)^{\bullet}$,
but their initial values differ.
By (\ref{eqn:Initial}) one has
\begin{eqnarray*}
&& \Phi_1(0, p)^{\bullet} = 1 + p_1 + 2p_1^2 + 3p_3 + 3p_1^3 + \cdots,
\end{eqnarray*}
hence by solving the systems of ODE's as in the genus zero case,
\begin{eqnarray*}
\Phi_1(\lambda, p)^{\bullet} & = & 1 +p_1 + 2p_2\sinh \lambda + 2p_1^2\cosh \lambda \\
&&+ p_3(4\cosh (3\lambda) -1) + 6p_1p_2\sinh (3\lambda) + p_1^3(2\cosh \lambda + 1) + \cdots.
\end{eqnarray*}
This matches with the results from \S \ref{sec:Example11}.

\subsection{Transfer to results on Hodge integrals}

As mentioned above,
one can transfer the above results to Hodge integrals by the ELSV formula.
This is rather straightforward so we will leave the detailed formulations to the reader.

\begin{example}
Taking $\eta_i = 1$ one get
\begin{eqnarray*}
&& \sum_{g \geq 0}
\frac{1}{l!} \prod_{i=1}^{l} \int_{\Mbar_{g, l}} \frac{\Lambda_g^{\vee}(1)}{\prod_{i=1}^{l} (1 - \psi_i)}
\cdot \lambda^{2g-2+2l} \\
& = & \sum_{ n \geq 1} \frac{(-1)^{n-1}}{n} \sum_{\eta}
\sum_{\substack{k_1 + \dots + k_n = l\\ k_1, \dots, k_n > 0}} \prod_{j=1}^n U_0^{(1^{k_j})}.
\end{eqnarray*}
Now
$$U_0^{(1^m)} = \sum_{|\rho| =m} \frac{\dim R_{\rho}}{m!} \exp\left[\begin{pmatrix} m\\2\end{pmatrix}
\frac{\chi_{\rho}(1^{m-2}2)}{\dim R_{\rho}}\lambda\right].$$
By changing $\exp$ to $\cosh$,
one recovers (4.36) in \cite{Mon-Son-Son}.
\end{example}

The results in \S \ref{sec:Example01} can be transferred to the following closed formulas for Hodge integrals:
\begin{eqnarray*}
&& \sum_{g \geq 0} \lambda^{2g} \int_{\Mbar_{g,1}} \frac{\Lambda_g^{\vee}(1)}{1-\psi_1} = 1, \\
&& \sum_{g \geq 0} \frac{2^2}{2!}
\int_{\Mbar_{g, l}} \frac{\Lambda^{\vee}_g(1)}{1 - 2 \psi_1}
\cdot \lambda^{2g-2+2 +1} = \frac{1}{2} \sinh \lambda, \\
&& \frac{1}{2} \sum_{g \geq 0}
\int_{\Mbar_{g, 2}} \frac{\Lambda_g^{\vee}(1)}{(1 - \psi_1)(1 - \psi_2)}
\cdot \lambda^{2g-2+2 +2}  =  \sinh^2 (\lambda/2), \\
&& \sum_{g \geq 0} \frac{3^3}{3!}
\int_{\Mbar_{g, l}} \frac{\Lambda^{\vee}_g(1)}{1 - 3 \psi_i}
\cdot \lambda^{2g-2+3 +1} = \frac{2}{9}\sinh^2(3\lambda/2)\\
&& \frac{2^2}{2!} \sum_{g \geq 0}
\int_{\Mbar_{g, 2}} \frac{\Lambda_g^{\vee}(1)}{(1 - \psi_1)(1 - 2\psi_2)}
\cdot \lambda^{2g-2+ 3 +2}  = \frac{1}{6}\sinh(3\lambda) - \frac{1}{2}\sinh \lambda, \\
&& \frac{1}{3!} \sum_{g \geq 0}
\int_{\Mbar_{g, 3}} \frac{\Lambda_g^{\vee}(1)}{\prod_{j=1}^3(1 - \psi_j)}
\cdot \lambda^{2g-2+ 3 +3}  =\frac{1}{9}\sinh^2(3\lambda/2) - \sinh^2(\lambda/2).
\end{eqnarray*}
By simple algebraic manipulations, they are equivalent to
\begin{eqnarray*}
&& \sum_{g \geq 0} \lambda^{2g} \int_{\Mbar_{g,1}} \frac{\Lambda_g^{\vee}(1)}{1-\psi_1} = 1, \\
&& \sum_{g \geq 0}
\int_{\Mbar_{g, l}} \frac{\Lambda^{\vee}_g(1)}{\frac{1}{2}(\frac{1}{2}- \psi_i)}
\cdot \lambda^{2g} = \cS(2\lambda), \\
&& \sum_{g \geq 0}
\int_{\Mbar_{g, 2}} \frac{\Lambda_g^{\vee}(1)}{(1 - \psi_1)(1 - \psi_2)}
\cdot \lambda^{2g}  = \frac{1}{2}\cS(\lambda)^2, \\
&& \sum_{g \geq 0}
\int_{\Mbar_{g, l}} \frac{\Lambda^{\vee}_g(1)}{\frac{1}{3}(\frac{1}{3}- \psi_1)}
\cdot \lambda^{2g} =\cS(3\lambda)^2\\
&& \sum_{g \geq 0}
\int_{\Mbar_{g, 2}} \frac{\Lambda_g^{\vee}(1)}{(1 - \psi_1)\cdot \frac{1}{2}(\frac{1}2-\psi_2)}
\cdot \lambda^{2g} = \frac{4}{3}\cS(2\lambda)^3 , \\
&& \sum_{g \geq 0}
\int_{\Mbar_{g, 3}} \frac{\Lambda_g^{\vee}(1)}{\prod_{j=1}^3(1 - \psi_j)}
\cdot \lambda^{2g}  =\frac{1}{2}(\cS(\lambda)^3\cS(3\lambda) + \cS(\lambda)^4),
\end{eqnarray*}
where
$$\cS(\lambda) = \frac{\sinh (\lambda/2)}{\lambda/2}.$$
Note the first four identities are known by other methods.
See e.g. \cite{Fab-Pan, Tia-Zho}.
The last identity matches with the identity
$$G_3(t, -1) = \frac{(2+ \cos t)}{3}\left(\frac{\sin (t/2)}{t/2}\right)^4$$
in \cite{Mon-Son-Son}.
The second to last identity seems to be new.
In higher degrees,
our method produces many more closed formulas for
$$\sum_{g \geq 0}
\int_{\Mbar_{g, l(\mu)}} \frac{\Lambda_g^{\vee}(1)}
{\prod_{j=1}^{l(\mu)} \frac{1}{\mu_j}(\frac{1}{\mu_j} - \psi_j)}
\cdot \lambda^{2g}.$$
We conjecture they all are polynomials in $\cS(\lambda), \dots, \cS(n\lambda), \dots$,
for all partitions of $d$.

{\bf Acknowledgements}.
{\em This version of the paper was prepared during the author's visit
to the Center of Mathematical Sciences, Zhejiang University.
The author wishes to thank the Center for the hospitality and the excellent
working environment.
He also likes to thank Professors Chiu-Chu Melissa Liu and Kefeng Liu
for discussions.
This work is partly supported by research grants from  NSFC and Tsinghua University}.

\end{document}